%

\magnification=\magstep1


\def\hexnumber#1{\ifcase#1 0\or1\or2\or3\or4\or5\or6\or7\or8\or9\or
	A\or B\or C\or D\or E\or F\fi }

\font\teneuf=eufm10
\font\seveneuf=eufm7
\font\fiveeuf=eufm5
\newfam\euffam
\textfont\euffam=\teneuf
\scriptfont\euffam=\seveneuf
\scriptscriptfont\euffam=\fiveeuf
\def\frak{\fam\euffam \teneuf}

\font\tenmsx=msam10
\font\sevenmsx=msam7
\font\fivemsx=msam5
\font\tenmsy=msbm10
\font\sevenmsy=msbm7
\font\fivemsy=msbm5
\newfam\msxfam
\newfam\msyfam
\textfont\msxfam=\tenmsx  \scriptfont\msxfam=\sevenmsx
  \scriptscriptfont\msxfam=\fivemsx
\textfont\msyfam=\tenmsy  \scriptfont\msyfam=\sevenmsy
  \scriptscriptfont\msyfam=\fivemsy
\edef\msx{\hexnumber\msxfam}

\mathchardef\upharpoonright="0\msx16

\def\Bbb#1{\tenmsy\fam\msyfam#1}

\def\Bigskip{\vskip2.2truecm}

\def\qed{{\vcenter{\hrule height.4pt \hbox{\vrule width.4pt height5pt
 \kern5pt \vrule width.4pt} \hrule height.4pt}}}
\def\notin{{\in}\kern-5.5pt / \kern1pt}
\def\ok{\vbox{\hrule height 8pt width 8pt depth -7.4pt
    \hbox{\vrule width 0.6pt height 7.4pt \kern 7.4pt \vrule width 0.6pt height 7.4pt}
    \hrule height 0.6pt width 8pt}}
\def\nt{{\leq}\kern-1.5pt \vrule height 6.5pt width.8pt depth-0.5pt \kern 1pt}
\def\sd{{\times}\kern-2pt \vrule height 5pt width.6pt depth0pt \kern1pt}
\def\zp#1{{\hochss Y}\kern-3pt$_{#1}$\kern-1pt}

\def\ZZ{{\Bbb Z}}
\def\sm{{\smallskip}}

\def\la{{\langle}}
\def\ra{{\rangle}}
\def\sub{\subseteq}

\def \o {\omega }

\font\capit=cmcsc10 scaled\magstep0

\font\bolds=cmssdc10 scaled\magstep0

\overfullrule=0pt
\openup1.5\jot

\def \o {\omega } \def \sub {\subseteq } \def \s {\sigma }

\def \l {\lambda }  
\def \s {\sigma } \def \m {\mu } \def \n {\nu } 
\def \g {\gamma }
\def \t {\tau } \def \a {\alpha } \def \b {\beta } \def \d {\delta } 
\def \z {\zeta }

\def \M {{\frak M}} \def \N {{\frak N}}
 
\font\gross=cmbx10 scaled \magstep1
\font\sgross=cmbx10 scaled \magstep2
  
\def \A {{\cal A}} \def \B {{\cal B}}

\noindent {\sgross The essentially free spectrum of a variety} 

\Bigskip

\noindent Alan H. Mekler\footnote{$^1$}{The author is supported by the
NSERC} 

\smallskip

\item{}{Department of Mathematics and Statistics,
Simon Fraser University, Burnaby, B.C. V5A 1S6 CANADA}

\smallskip

\noindent Saharon Shelah\footnote{$^2$}{The author is supported by the
United States--Israel Binational Science Foundation;
publication 417.}

\smallskip

\item{}{Institute of Mathematics, Hebrew University, Givat Ram, 91904
Jerusalem, ISRAEL}

\smallskip

\noindent Otmar Spinas\footnote{$^3$}{The author is supported by the
Schweizer Nationalfonds.}

\smallskip



\item{}{Department of Mathematics, University of California, Irvine,
CA 92717, USA

\Bigskip

{\narrower{ABSTRACT:  We partially prove a conjecture from [MeSh]
which says that the spectrum of almost free, essentially free,
non-free algebras in a variety is either empty or consists of the
class of all successor cardinals.

}}

\Bigskip

{\gross Introduction and notation}

\Bigskip

Suppose that $T$ is a {\bf variety} in a countable vocabulary $\t $. This
means that $\t $ is a countable set of function symbols and $T$ is a set of
equations, i.e. sentences of the form $\forall x_1,\dots ,x_n$ $(\s
_1(x_1,\dots ,x_n)= \s _2(x_1,\dots ,x_n))$ where $\s _i$ are $\t
-$terms. The class of all models of $T$ will be denoted by Mod$(T)$,
and a member of Mod$(T)$ is called an {\bf algebra} in the variety $T$. Let
$M\in $ Mod$(T)$. For $A\sub M$, $\langle A \ra $ denotes the submodel
of $M$ generated by $A$. Such $A$ is called a {\bf free basis} (of $\la A\ra
$) if no distinct $a_1,\dots ,a_n\in A$ satisfy an equation which is
not provable from $T$. Moreover, $M$ is called {\bf free} if there exists a
free basis of $M$, i.e. one which generates $M$. By $F_\l $ we denote
the free algebra with free basis of size $\lambda $, where $\l $ is a
cardinal. For $M_1,M_2\in $ Mod$(T)$, the {\bf free product} of $M_1$ and
$M_2$ is denoted by $M_1\ast M_2$. Formally it is obtained by building
all formal terms in the language $\t $ with constants belonging to the
disjoint union of $M_1$ and $M_2$, and then identifying them according
to the laws in $T$. For $M,$ $\la M_\nu : \nu < \a \ra $ such that $M,
M_\nu \in $ Mod$(T)$ and $M$ is a submodel of $M_\nu $ for all $\nu
<\a $, the {\bf free product of the} $M_\nu '$s {\bf over} $M$ is defined
similarly, and it is denoted by $\ast _M\{ M_\nu :\nu <\a \} $; the
intention being that distinct $M_\nu $,
$M_{\nu '} $ are disjoint outside $M$ except for those equalities
which follow from
the laws in $T$ and the equations in Diag$(M_\nu )\cup $
Diag$(M_{\nu '})$. Here Diag denotes the {\bf diagram} of a model.
For $M,N\in $ Mod$(T)$ we say ``$N/M$ is free'' if
$M$ is a submodel of $N$ and there exists a free basis $A$ {\bf of}
$N$ {\bf over} $M$, i.e. $A$ is a free basis, $N=\la M\cup A\ra $ and
between members of $\la A\ra $ and $M$ only those equations hold which
follow from $T$ and Diag$(M)$.

Suppose $|M|=\l $. Then $M$ is called {\bf almost free} if there
exists an increasing continuous family $\la M_\nu :\nu < \hbox{ cf}(\l
)\ra $ of free submodels of size $<\l $ with union $M$. Moreover, $M$
is called {\bf essentially free} if there exists a free $M'\in $
Mod$(T)$ such that $M\ast M'$ is free, {\bf essentially non-free}
otherwise. The {essentially free spectrum} of the variety $T$ which is
denoted by EINC$(T)$, is the class of cardinals $\l $ such that there
exists $M\in $ Mod$(T)$ of size $\l $ which is almost free and
essentially free, but not free. 

In [MeSh] the {\bf essentially non-free spectrum}, i.e. the spectrum of 
cardinalities
of almost free and essentially non-free algebras in a variety $T$, has been
investigated, and it is shown that this spectrum has no simple
description in ZFC, in general. Here we will show that the situation
is different for EINC($T$). Firstly, by a general compactness theorem
due to the second author (see [Sh]), EINC($T$) contains only regular
cardinals. Secondly, we will show that EINC$(T)$ is contained in the
class of successor cardinals. Our conjecture is that EINC$(T)$ is
either empty or equals the class of all successor cardinals (depending
on $T$). Motivating examples for this conjecture are among others 
$\ZZ /4\ZZ-$modules (where EINC is empty) and $\ZZ /6\ZZ -$modules (where EINC
consists of all successor cardinals) (see [EkMe, p.90]).
We succeed to prove the conjecture to a certain extent.
Namely, we prove the following theorem.

\sm

{\bolds Theorem.}  {\it If for some cardinal $\mu $, 
$(\mu  ^{\aleph _0})^+ \in $} EINC$(T)$, {\it then every successor 
cardinal belongs
to} EINC$(T)$.

\sm

For the proof we will isolate a property of $T$, denoted Pr$_1(T)$,
which says that a countable model of $T$ with certain properties
exists, and then show that, on the one hand, the existence
of $M\in $ Mod$(T)$ in any cardinality of the form $(\mu ^{\aleph _0})^+$
implies that Pr$_1(T)$ holds, and on the other hand, from Pr$_1(T)$ an
algebra $M\in $ Mod$(T)$ can be constructed in every successor cardinality.

\Bigskip 

{\gross 1. EINC({\it T}) is contained in the class of successor
cardinals}

\Bigskip

{\bolds Theorem.}  {\it For every variety $T$, EINC$(T)$ is contained
in the class of successor cardinals.}

\sm

{\capit Proof:}  Suppose $\l \in $ EINC$(T)$. By the main result of
[Sh], $\l $ must be regular. So suppose $\l $ is a regular limit
cardinal. Let $\M \in $ Mod$(T)$ be generated by $\{ a_\a : \a < \l \}
$ and suppose that $\M $ is almost free and essentially free. We will
show that then $\M $ must be free, and hence does not exemplify $\l
\in $ EINC$(T)$. 

By assumption and a L\"owenheim-Skolem argument, $\M \ast F_\l $ is
free. Let $\{ c_\nu :\nu <\l \} $, $\{ b_\nu :\nu <\l \} $ be a free basis
of $\M \ast F_\l $, $F_\l $, respectively.

Let $\chi $ be a large enough regular cardinal, and let $C\sub \l $ be the
club
consisting of all $\a $ such that for some substructure ${\cal A}\prec
\la H(\chi ), \in , \prec _\chi \ra $ of size $<\l $ which contains
$\M , F_\l , \{ a_\nu :\nu <\l \} , \{ b_\nu : \nu <\l \} $ and $\{
c_\nu : \nu <\l \} $, we have ${\cal A}\cap \l =\a $. Here $H(\chi )$ is the
set of all sets which are hereditarily of cardinality $<\chi $, and 
$\prec _\chi $ is a fixed well-ordering of $H(\chi )$. 
Note that the
information about $\M $ reflects to each $\a \in C$, especially $\la
\{ c_\nu :\nu <\a \} \ra =\la \{ a_\nu :\nu <\a \} \ra \ast \la \{
b_\nu :\nu <\a \} \ra $.

Since $\M $ is supposed to be almost free, the set 

$$C_0=\{ \a \in C: \la \{ a_\nu :\nu <\a \} \ra \hbox{ is free} \} $$

\noindent is still a club. Let $\a ,\b \in C_0$ be cardinals with $\a
<\b $. We will show that $\la \{ a_\nu :\nu <\b \} \ra / \la \{ a
_\nu : \nu <\a \} \ra $ is free. This will suffice to conclude that
$\M $ is free since the cardinals below $\l $ are a club and hence 
$C_1 =\{ \a \in C_0: \a \hbox{ is a cardinal} \} $ is a club such that
for every $\a ,\b \in  C_1$ with $\a <\b $, $\la \{ a_\nu :\nu <\b \}
\ra / \la  \{ a_\nu : \nu <\a \} \ra $ is free. 

For the proof, let $\{ d_\nu : \nu < \b \} $ be a free basis of $\la
\{ a_\nu : \nu < \b \} \ra $. As $\a =|\a |<|\b |=\b $ we may assume
$\la \{ a_\nu :\nu <\a \} \ra \sub \la \{ d_\nu : \nu <\a \} \ra $.
Hence easily

$$\la \{ a_\nu : \nu <\b \} \ra \cong _{\la \{ a_\nu :\nu <\a \} \ra }
\la \{ a_\nu :\nu <\b \} \ra \ast F_\b \, ,$$

\noindent i.e. there exists an isomorphism which leaves $\la \{ a_\nu
:\nu <\a \} \ra $ fixed. But $\la \{ a_\nu :\nu <\b \} \ra \ast F_\b
\cong \la \{ c_\nu :\nu <\b \} \ra $ and $\la \{ c_\nu :\nu <\b \} \ra
/ \la \{ c_\nu :\nu <\a \} \ra $ is free. Moreover $\la \{ c_\nu :\nu
<\a \} \ra =\la \{ a_\nu :\nu <\a \} \ra \ast \la \{ b_\nu :\nu < \a
\} \ra $ and hence $\la \{ c_\nu :\nu <\a \} \ra / \la \{ a_\nu :\nu
<\a \}  \ra $ is free. Consequently $\la \{ a_\nu :\nu <\b \} \ra /
\la \{ a_\nu :\nu <\a \} \ra $ is free. $\qed $ 

\Bigskip

{\gross 2.  EINC({\it T}) is either empty or contains almost all
successor cardinals}

\Bigskip
  
{\bolds Definition 2.1.}  The property Pr$_1(T)$ says: There exist
$N,M\in $ Mod$(T)$ such that $N$ is countably generated, $M$ is a subalgebra of
$N$ and the following 
clauses hold:

\itemitem{$(i)$} $M$ has a free basis;

\itemitem{$(ii)$} $N\ast F_{\aleph _0}$/$M$ is free;

\itemitem{$(iii)$} $\ast _M\{N:n\in \o \} \ast F_{\aleph _0}$/$M\ast
F_{\aleph _0}$ is not free.

\sm

{\bolds Theorem 2.2.}  {\it Suppose that} Pr$_1(T)$ {\it holds and $\l
$ is a successor cardinal. Then} $\l \in $ EINC$(T)$.

\sm

{\capit Proof:}  Let $\l =\mu ^+$. Let $M,N$ witness Pr$_1(T)$. Let
$\N =  \ast _M \{
N:\a < \l \} $. We claim that $\M = \N \ast F_\m $ exemplifies that
$\l \in $ EINC$(T)$. Let $\{ c_\a :\a <\m \} $ be a free basis of
$F_\m $.

Firstly, $\M $ is almost free: For $\a <\l $ let $\N _\a = \ast _M \{
N:\n <\a \}$. Then clearly $\la \N _\a \ast F_\m :\a <\l \ra $ is a
$\l -$filtration of $\M $. Moreover $\N _\a \ast F_\m $ is free for
every $\a <\l $, since easily $\N _\a \ast F_\m \cong \ast _M \{ N\ast
F_{\aleph _0}:\n <\a \} $ and by Pr$_1(T)$, $M$ is free and $N\ast
F_{\aleph 0}$/$M$ is free.

Secondly, $\M \ast F_\l \cong \N \ast F_\l $ is free, since $\N \ast
F_\l \cong \ast _M \{ N\ast F_{\aleph _0}:\a <\l \} $ is free as in
the proof of almost freeness.

Thirdly, $\M $ is not free. By contradiction, suppose that $I=\{ d_\n
:\n <\l \} $ were a free basis of $\M $. 

Let $\chi $ be a large enough regular cardinal, and let $\A \prec \la H(\chi
),\in ,\prec _\chi \ra $ such that $|\A |=\m $, $\m +1\sub \A $, and
$\l ^+, N,M, \N ,\M ,F_\m ,I\in \A $. Next choose $\B \prec \la H(\chi
),\in ,\prec _\chi \ra $ such that $|\B |=\aleph _0 $, and
$\A ,\l ^+, N,M, \N ,\M ,F_\m ,I\in \B $. 

Let $u=\B \cap \l \setminus (\A \cap \l )$, $v=\A \cap \B \cap \l $,
$w=\A \cap \B \cap \m $. Notice that $w=\B \cap \m $. Define $M_1 =\A
\cap \B \cap \M $. Now easily $M_1$ is countably generated and it has
the form

$$M_1= \ast _M \{ N:\a \in v \} \ast \la \{ c_\a :\a \in w \} \ra .$$

\noindent Hence $M_1 \cong _M \ast _M \{ N:n\in \o \} \ast F_{\aleph
_0} \cong _M \ast _M \{ N\ast F_{\aleph _0} :n\in \o \} \cong _M M\ast
F_{\aleph _0}$, where for the last isomorphy we applied $(ii)$ from
Pr$_1(T)$. Next define $M_2= \B \cap \M $. Then easily

$$M_2 =\ast _M \{ N:\a \in u \} \ast _M M_1 .$$

\noindent Hence by the isomorphy above we have

$$M_2\cong \ast _M \{ N:n\in \o \} \ast F_{\aleph _0} .$$

By $(iii)$ from Pr$_1(T)$ we conclude that $M_2$/$M_1$ is not free. On
the other hand, $\{ d_\n :\n \in u \} $ witnesses  that $M_2$/$M_1$ is
free, a contradiction. $\qed $

\sm

{\bolds Theorem 2.3.}  {\it Suppose $\l , \m $ are cardinals such that
$\l = \m ^+$, $\m ^{\aleph _0} =\m $ and} $\l \in $ EINC$(T)$. {\it Then
Pr$_1(T)$ holds.}

\sm

{\capit Proof:}  Let ${\M }$ exemplify $\l \in $ EINC$(T)$. Let $\{
a_\n :\n < \l \} $ generate $\M $. Let $F$ be free such that $\M \ast
F$ is free. Without loss of generality we may assume that $F=F_\l $;
in fact, if $|F|<\l $ then we may replace $F$ by $F\ast F_\l $ which
is isomorphic to $F_\l $, and if $|F|>\l $ use a L\"owenheim--Skolem
argument. So let $\{ b_\n :\n < \l \} $ be a free basis of $F$, and
let $\{ c_\nu :\n < \l \} $ be a free basis of $\N =\M \ast F$. 

Let $\chi $ be a large enough regular cardinal, and let $N_\a $, for every $\a
< \l $, be a countable elementary substructure of $\la H(\chi ), \in ,
\prec _\chi \ra $ such that $\a , \M ,F,\N , \{ a_\n :\n <\l \} , \{
b_\n :\n <\l \} , \{ c_\n :\n <\l \} $ belong to $N_\a $. Let $u_\a
=N_\a \cap \l $. 

By assumption on $\M $ ($\M $ is almost free), the set

$$\vbox{\halign{\hfill #&#\hfill \cr

$\{ \a < \l :$ &$ \la \{ a_\n :\n <\a \} \ra \hbox{ is free } \wedge
\la \{ c_\n :\n < \a \} \ra  = \la \{ a_\n :\n <\a \} \ra  \ast
\la \{ b_\n :\n <\a \} \ra \} $\cr }}$$

\noindent contains a club; let $C$ be the $\prec _\chi -$least one.
Hence $C\in N_\a $ for every $\a <\l $.

Using elementarity, it is easy to see that for every $\a \in C$ the
following three clauses hold ((1) holds for every $\a <\l $):

\item{(1)}{$\la \{ c_\n :\n \in u_\a \} \ra =\la \{ a_\n :\n \in u_\a
\} \ra \ast \la \{ b_\n :\n \in u_\a \} \ra $;}

\item{(2)}{$\la \{ c_\n :\n \in u_\a \cap \a \} \ra =\la \{ a_\n :\n
\in u_\a \cap \a
\} \ra \ast \la \{ b_\n :\n \in u_\a \cap \a \} \ra $;} 

\item{(3)}{$\{ \la a_\n :\n \in u_\a \cap \a \} \ra $ is free, and 
$\la \{ a_\n :\n \in \a \} \ra \, $/$\,\{ \la a_\n :\n \in u_\a \cap
\a \} \ra $ is free.}

\sm
To prove (3), let $\la d_\n :\n \in I\ra $ be the $\prec _\chi -$least
free basis of $\la \{ a_\n :\n \in \a \} \ra $. So by elementarity
$\la d_\n : \n \in I\ra \in N_\a $ and $\la \{ d_\n :\n \in I\cap N_\a
\} \ra =\la \{ a_\n :\n \in u_\a \cap \a \} \ra $. Hence $\{
d_\n :\n \in I\cap N_\a \} $ and $\{ d_\n : \n \in I\setminus N_\a \}
$ witness that (3) holds.

Moreover it is not difficult to see that $C_0= \{ \a \in C: \a =
\bigcup \{ u_\n :\n <\a \} \} $ is still a club. Hence $S_0 =\{ \a \in
C_0 : \hbox{ cf}(\a )>\o \} $ is stationary. By Fodor's Lemma, for
some $\a ^*<\l $, $S_1 =\{ \a \in S_0: u_\a \cap \a \sub \a ^* \} $ is
stationary. By assumption, $|\a ^*|^{\aleph _0} \leq \m ^{\aleph _0}
<\l $. So by thinning out $S_1$ further (using this assumption and the
$\l -$completeness of the nonstationary ideal on $\l $), we may find a
stationary $S_2\sub S_1 $ and $u^*\sub \a ^*$ such that for every $\d
_1 ,\d _2\in S_2$ the following hold:

\sm
\item{(4)}{$u_{\d _1}\cap \d _1 =u^*$;}

\item{(5)}{o.t.$(u_{\d _1} )=$ o.t.$(u_{\d _2} )$, and the unique
order-preserving map $h=h_{\d _1 \d _2} : u_{\d _1}\rightarrow u_{\d
_2} $ induces (by $c_\n \rightarrow c_{h(\n )}$) an isomorphism from $\la \{
c_\n :\n \in u_{\d _1} \} \ra $ onto $ \la \{ c_\n :\n \in u_{\d _2}
\} \ra $ which maps $a_\n $ to $a_{h(\n )}$ and $b_\n $ to $b_{h(\n )}$.}  

\sm

Let $\d ^*=\min (S_2 \setminus \m )$, $M=\la \{ a_\n :\n \in u^* \}
\ra $ and $N=\la \{ a_\n :\n \in u_{\d ^*} \} \ra $.

As $\d ^* \in C$, by elementarity we know that $M$ is free.

As $\{ c_\n :\n \in \l \} $ is a free basis, clearly $\la \{ c_\n :\n
\in u_{\d ^*} \} \ra \, $/$\, \la \{ c_\n : \n \in u^* \} \ra $ is
free, and by (2) and as $\d ^* \in S_2\sub C$, also $\la \{ c_\n : \n \in u^*
\} \ra \, $/$\, M$ is free. Finally, $\la \{ c_\n :\n \in u_{\d ^*}\}
\ra  \cong N\ast F_{\aleph _0}$ by (1). Hence we conclude that $N\ast F_{\aleph
_0} \, $/$\, M$ is free. 

Hence, if the pair $M,N$ does not exemplify Pr$_1(T)$, then $(iii)$ in
its definition fails. We will use this to show that then $\M $ is
free, which contradicts our assumption. Then we conclude that
Pr$_1(T)$ holds.

By induction on $\z < \l $ we choose $w_\z \sub \l $ such that the
following requirements are satisfied:

\item{(6)}{$w_0=\d ^*$;}

\item{(7)}{$|w_\z |< \l $;}

\item{(8)}{for $\z $ limit, $w_\z =\bigcup \{ w_\n :\n < \z \} $;}

\item{(9)}{if $\gamma (\z )=\min (\l \setminus w_\z )$, then $w_{\z
+1}=w_\z \cup \{ \g (\z )\} \cup \{ \b (\z ,n):n\in \o \} $, where the
$\b (\z ,n)$ belong to $S_2$, and for any $m,n\in \o $ with $m<n$,
$\bigcup \{ u_{\g (\z )}, u_\n : \n \in w_\z \} < \min (u_{\b (\z ,n)}
\setminus u^*)$ and $\sup (u_{\b (\z ,m)})<\min (u_{\b (\z
,n)}\setminus u^*)$ hold.}

By (6) and $\d ^*\in C_0\sub C$ we conclude that $\bigcup \{ u_\n :\n
\in w_0\} =\d ^* $ and $\la \{ a_\a :\a \in \d ^* \} \ra $ is free. By
(8) and (9) it is clear that the sequence

$$\la \la \{ a_\a :\a \in \bigcup \{u_\n : \n \in w_\z \} \} \ra : \z
< \l \ra $$

\noindent is increasing and continuous with limit $\M $. Hence the
following claim gives the desired contradiction:

\sm

{\bolds Claim.}  {\it For every} $\z < \l $, $\la \{ a_\a :\a \in
\bigcup  \{u_\n : \n \in w_{\z +1} \} \} \ra $/$\la \{ a_\a :\a
\in \bigcup \{u_\n : \n \in w_\z \} \} \ra  $ {\it is free.}

\sm

{\capit Proof:}  Let us introduce the following notation. For $x\in \{
a,b,c \} $ and $I\sub \l $ set:

$$\vbox{\halign{\hfill #&#\hfill \cr
$Z^x_I$&$=\la \{ x_\a : \a \in \bigcup \{ u_\n :\n \in I\} \} \ra $\cr
$W^x_\z $&$=Z^x_{w_\z }$\cr
$K^x$&$=\la \{ x_\a :\a \in u^* \} \ra $, so $K^a=M$.\cr }}$$

The Claim will follow from the following three facts:

\item{(10)}{$Z^c_I=Z^a_I \ast Z^b_I$;}

\item{(11)}{$\la W^a_\z \cup Z^a_{\{ \g (\z )\} } \ra \ast F_{\aleph
_0}$/$W^a_\z $ is free;}

\item{(12)}{$W^a_{\z +1}=\ast _{K^a} \{  Z^a_{w_\zeta \cup \{ \g (\z
)\} }  , Z^a_{\{ \b (\z ,n)\} }: n\in \o \} $.}

For (10), to prove $Z^c_I= \la Z^a_I\cup Z^b_I\ra $ is rather
straightforward by using (1).
Moreover there exists a homomorphism $h: \M \ast F \rightarrow Z^c_I $
which is the identity on $Z^c_I$ and maps $\M $ onto $Z^a_I$; $h$ can be
defined by letting

$$h(c_\a )=\cases{c_\a & if $\a \in \bigcup \{ u_\n :\n \in I\} $,\cr
                  a_0& otherwise.\cr} $$

Now suppose that $\la Z^a_I\cup Z^b_I\ra \models \phi (\bar a,\bar
b)$, where $\phi $ is an equation and $\bar a \sub \{ a_\a :\a \in
\bigcup \{ u_\n :\n \in I\} \} $, $\bar b \sub \{ b_\a :\a \in 
\bigcup \{ u_\n :\n \in I\} \} $ are finite. Then this equation holds
in $\M \ast F$, of course. As $\{ b_\n :\n <\l \} $ is a free basis of
$F$, we conclude that $\phi (\bar a,\bar b)$ is provable from finitely
many equations in Diag$(\M )$ and the laws of the variety. But $h$
maps this proof to a proof from Diag$(Z^a_I)$ and leaves $\bar a,\bar
b$ fixed. Consequently $Z^c_I=Z^a_I\ast Z^b_I$ holds.

To prove (11), first clearly $Z^c_{w_\z \cup \{ \g (\z )\}}$/$W^c_\z $
is free. Hence by (10), $Z^c_{w_\z \cup \{ \g (\z )\} }$/$W^a_\z $ is
free. We may assume that $|w_\z |=\m $; hence $Z^a_{w_\z \cup \{ \g
(\z )\} }\ast F_\m $/$W^a_\z $ is free. How to get $\m $ down to
$\aleph _0$? Choose $\bar N\prec (H(\chi ),\in ,\prec _\chi )$ such that
$\bar N$ is countable and contains everything relevant, especially $W^a_\z
$, $Z^a_{\g (\z )}$, $F_\m $. Let $X\in \bar N$ be a free basis of
$Z^a_{w_\z \cup \{ \g (\z )\} } \ast F_\m $ over $W^a_\z $. Then
$X\cap \bar N$ is a free basis of $\bar N\cap (Z^a_{w_\z \cup \{ \g (\z )\} }
\ast  F_\m )$ over $\bar N\cap W^a_\z $. Moreover $\bar N\cap
(Z^a_{w_\z \cup \{ 
\g   (\z )\} } \ast F_\m ) = \la (\bar N\cap W^a_\z )\cup Z^a_{\{ \g
(\z )\}  }\ra
\ast (\bar N\cap F_\m )\cong \la (\bar N\cap W^a_\z )  \cup Z^a_{\{
\g  (\z )\} }\ra
\ast F_{\aleph _0} $. We claim that $X\cap \bar N$ is a free basis over
$W^a_\z $ of what it generates over $W^a_\z $, namely 
$Z^a_{w_\z \cup \{ \g (\z )\} }
\ast  F_{\aleph _0} $. Otherwise there were finite sets $X_0\sub X$
and $Y_0\sub W^a_\z $ such that $X_0\cup Y_0$ satisfies an equation
which does not follow from the laws of the variety and the equalities
in Diag$(W^a_\z )$. By elementarity we can find $Y_1\sub W^a_\z \cap
\bar N$ such that $X_0\cup Y_1$ satisfies the same equation, a
contradiction.

To prove (12), if $a$ is replaced by $c$, then (12) is easily verified
by using the free basis $\{ c _\n : \n <\l \} $. But then using (10) and
$K^c=K^a\ast K^b$ we easily finish.

Finally, as $\d ^* \in C$ and $\d ^*$ has uncountable cofinality, by
(3) we may choose $F_{\aleph _0}\sub \la \{ a_\n :\n \in \d ^* \} \ra
$ (i.e. an algebra isomorphic to $F_{\aleph _0}$) such that $M\cap
F_{\aleph _0} =\emptyset $ and even $\la M\cup F_{\aleph _0}\ra =
F_{\aleph _0}\ast M $. By (12) we conclude

$$W^a_{\z +1}= Z^a_{w_\z \cup \{ \g (\z )\} }\ast _{F_{\aleph
_0}\ast M }  (F_{\aleph _0} \ast  \ast _M\{Z^a_{\{ \b (\z ,n)\} }:n\in 
\o \}).$$

Moreover by construction (as $\b (\z ,n)\in S_2$), $Z^a_{\{ \b (\z
,n )\} }\cong  _M N$ for every $n\in
\o $. Hence $F_{\aleph _0} \ast  \ast _M\{Z^a_{\{ \b (\z ,n)\} }:n\in \o \}
\cong _M F_{\aleph _0} \ast  \ast _M\{ N:n\in \o \} $. By assumption
$F_{\aleph _0} \ast  \ast _M\{ N:n\in \o \} $/$F_{\aleph _0}\ast M$ is free,
and so clearly of rank $\aleph _0$. We conclude that $W^a_{\z +1}
\cong _{W^a_\z } Z^a_{w_\z \cup \{ \g (\z )\} }\ast  F_{\aleph _0}$, and so
$W^a_{\z +1}$/$W^a_\z $ is free by (11). $\qed $ $\qed $

\Bigskip

\centerline{{\gross References}}

\bigskip

\item{[EkMe]}{P.C. Eklof and A.H. Mekler, Almost free modules:
set-theoretic methods, North-Holland 1990}

\smallskip

\item{[MeSh]}{A.H. Mekler and S. Shelah, Almost free algebras, Israel
J. Math., to appear}

\smallskip

\item{[Sh]}{S. Shelah, A compactness theorem for singular cardinals,
free algebras, Whitehead problem and transversals, Israel J. Math.
21(1975), 319--349}

\bye